\documentclass[runningheads,a4paper]{llncs}

\usepackage{amssymb}
\usepackage{amsmath}
\usepackage{graphicx}
\newcommand{\tensor}[1]{\boldsymbol{\mathcal{#1}}}

\newcommand{\mat}[1]{\mathbf{#1}}
\newcommand{\vect}[1]{\mathbf{#1}}

\usepackage{url}

\newcommand{\keywords}[1]{\par\addvspace\baselineskip
\noindent\keywordname\enspace\ignorespaces#1}

\usepackage{makecell}
\usepackage{float}

\begin{document}

\mainmatter 

\title{Completion of High Order Tensor Data with Missing Entries via Tensor-train Decomposition}

\titlerunning{Completion of High Order Tensor Data with Missing Entries}

\author{Longhao Yuan$^{1,2}$ \and Qibin Zhao$^{2,3}$ \and Jianting Cao$^1$}
\institute{Graduate School of Engineering, Saitama Institute of Technology, Japan
\and Tensor Learning Unit, RIKEN Center\\ for Advanced Intelligence Project (AIP), Japan
\and School of Automation, Guangdong University of Technology, China\\ \url{{longhao.yuan, qibin.zhao}@riken.jp, cao@sit.ac.jp}}

\authorrunning{Longhao Yuan, Qibin Zhao and Jianting Cao}

\maketitle

\begin{abstract}
In this paper, we aim at the completion problem of high order tensor data with missing entries. The existing tensor factorization and completion methods suffer from the curse of dimensionality when the order of tensor $N>>3$.  To overcome this problem, we propose an efficient algorithm called TT-WOPT (Tensor-train Weighted OPTimization) to find the latent core tensors of tensor data and recover the missing entries. Tensor-train decomposition, which has the powerful representation ability with linear scalability to tensor order, is employed in our algorithm. The experimental results on synthetic data and natural image  completion demonstrate that our method significantly outperforms the other related methods. Especially when the missing rate of data is very high, e.g., 85\% to 99\%, our algorithm can achieve much better performance than other state-of-the-art algorithms.

\end{abstract}

\section{Introduction}

Tensor is a high order generalization of vectors and matrices, which is suitable for natural data with the characteristic of multi-dimensionality. For example, a RGB image can be represented as a three-way tensor: $height \times width \times channel$ and a video sequence can be represented by a $height \times width \times channel \times time$ form data. When the original data is transformed into matrix or vector forms, the structure information and adjacent relation of data will be lost. Tensor is the natural representation of data that can retain the high dimensional structure of data. In recent decades, tensor methodologies have attracted a lot of interests and have been applied to various fields such as image and video completion \cite{jour1,jour2}, signal processing \cite{jour3,jour4}, brain computer interface \cite{jour5}, image classification \cite{Proceeding1,Proceeding2} , etc. Many theories, algorithms and applications of tensor methods have been proposed and studied, which can be referred in the comprehensive review \cite{jour6}. 

Most tensor decomposition methods assume that the tensor has no missing entries and is complete. However, in practical situations, we may encounter some transmission or device problems which result in that the collected data has missing and unknown entries. To solve this problem, the study on high order tensor decomposition/factorization with missing entries becomes significant and has a promising application aspect. The goal of tensor decomposition of missing data is to find the latent factors of the observed tensor, which can thus be used to reasonably predict the missing entries. The two most popular tensor decomposition methods in recent years are CANDECOMP/PARAFAC(CP) decomposition \cite{jour7,jour8} and Tucker decomposition \cite{jour9}. There are many proposed methods that use CP decomposition to complete data with missing entries. CP weighted optimization (CP-WOPT) \cite{jour1} applies optimization method to finding the optimal CP factor matrices from the observed data.  Bayesian CP factorization \cite{jour2} exploits Bayesian probabilistic model to automatically determine the rank of CP tensor while finding the best factor matrices. The method in \cite{jour10} recovers low-n-rank tensor data with its convex relaxation by alternating direction method of multipliers (ADM). 

However, because of the peculiarity of CP and Tucker model, they can only reach a relatively high accuracy in low-dimension tensors. When it comes to a very high dimension, the performance of applying these models to missing data completion will decrease rapidly. As mentioned above, many natural data's original form is high dimension tensor, so the models which are not sensitive to dimensionality should be applied to perform the tensor decomposition. In this paper, we use tensor-train decomposition \cite{jour11} which is free from the curse of high dimension to perform tensor data completion. Our works in this paper are as follows: (a) We develop a optimization algorithm named tensor-train weighted optimization (TT-WOPT) to find the factor core tensors of tensor-train decomposition. (b) By TT-WOPT algorithm, tensor-train decomposition model is applied to incomplete tensor data. Then the factor core tensors are calculated and used to predict the missing entries of the original data. (c) We conduct simulation experiments to verify the accuracy of our algorithm and compare it to other algorithms. In addition, we carry out several real world experiments by applying our algorithm and other state-of-the-art algorithms to a set of $256 \times 256 \times 3$ images with missing entries. The experiment results show that our method performs better in image inpainting than other state-of-the-art approaches. In addition, by converting the image of size $256 \times 256 \times 3$ to a much higher dimension, our algorithm can successfully recover images with 99\% missing entries while other existing algorithms fail at this missing rate. These results demonstrate that tensor-train decomposition with high order tensorizations can achieve high compressive and representation abilities.

\section{Notations and Tensor-train Decomposition}

\subsection{Notations}
In this paper, vectors are denoted by boldface lowercase letters, e.g., $\vect{x}$. Matrices are denoted by boldface capital letters, e.g., $\mat{X}$. Tensors of order $N\geq 3$ are denoted by Euler script letters, e.g., $\tensor{X}$. $\mat{X}^{(n)}$ denotes the $n$th matrix of a matrix sequence and the representation of vector and tensor sequence is denoted by the same way. When the tensor $\tensor{X}$ is in the space of $\tensor{X}  \in\mathbb{R}^{I_1\times I_2\times\cdots \times I_N}$, $\mat{X}_{(n)}$ denotes the n-mode matricization of $\tensor{X}$, see \cite{jour6}. The $(i_{1},i_{2},\cdots,i_{N})$th element of $\tensor{X}$ is denoted by $x_{i_{1}i_{2}\cdots i_{N}}$ or $\tensor{X}(i_{1},i_{2},\cdots,i_{N})$.\\

\subsection{Tensor-train Decomposition}
The most important feature of tensor-train decomposition is that no matter how high the dimension of a tensor is, it decomposes the tensor into a sequence of three-way tensors. This is a great advantage in modeling high dimension tensor because the number of model parameters will not grow exponentially by the increase of the tensor dimension. For example, the number of parameters in Tucker model is $\tensor{O}(NIR+R^{N})$ where $N$ is the number of dimension, $R$ is the size of Tucker core tensor and $I$ is the size of each dimension of the tensor. For tensor-train decomposition, the number of parameters is $\tensor{O}(NIr^2)$ where $r$ is rank of TT-tensor. Therefore, TT-model needs much fewer model parameters than Tucker model.

Tensor-train decomposition is to decompose a tensor into a sequence of tensor cores. All the tensor cores are three-way tensors. In particular, the TT decomposition of a tensor $\tensor{X}  \in\mathbb{R}^{I_1\times I_2\times\cdots \times I_N}$ is expressed as follow:
\begin{equation}
\label{tt_decom}
\tensor{X}=\ll \tensor{G}^{(1)},\tensor{G}^{(2)},\cdots,\tensor{G}^{(N)} \gg,
\end{equation}
where  $\tensor{G}^{(1)},\tensor{G}^{(2)},\cdots,\tensor{G}^{(N)} $ is a sequence of three-way tensor cores with size of $1\times I_{1} \times r_{1},r_{1} \times I_{2} \times r_{2}, \cdots , r_{N-1} \times I_{N} \times 1$. The sequence $\{1, r_{1},r_{2},\cdots ,r_{N-1},1\}$ is named TT-ranks which can limit the size of every core tensor. Each element of tensor $\tensor{X}$ can be written as the following index form:
\begin{equation}
\label{TT_index}
\ x_{i_{1}i_{2}\cdots i_{N}}=\mat{G}^{(1)}_{i_{1}}\times \mat{G}^{(2)}_{i_{2}} \times \cdots \times \mat{G}^{(N)}_{i_{N}},
\end{equation} 
where $\mat{G}^{(n)}_{i_{n}}$ is the $i_{n}$th slice of the $n$th core tensor. See the concept of slice in \cite{jour6}. 

Currently, there is few study about how to compute TT-ranks efficiently. In paper \cite{jour11} where tensor-train decomposition is proposed, the author advances an algorithm named TT-SVD to calculate the core tensors and TT-ranks. Although it has the advantage of high accuracy and high efficiency, the TT-ranks in the middle core tensors must be very high to compensate the low TT-ranks in the border core tensors, which leads to the unreasonable distribution of TT-ranks and redundant model parameters. Therefore, the TT-ranks calculated by TT-SVD may not be the optimal one. In this paper, we manually set the TT-ranks to a smooth distribution and use TT-WOPT algorithm to calculate the core tensors. Though we do not have a good TT-rank choosing strategy, much fewer model parameters are needed. The simulation results and experiment results also show high accuracy and performance.

\section{TT-WOPT Algorithm}
Most of the tensor decomposition methods, which are used for finding the latent factors, only aim at the fully observed data. When data has missing entries, we cannot use these methods to predict the missing entries. Weighted optimization method minimizes the distance between weighted real data and weighted optimization objective. When the optimization is finished, it means the obtained tensor decomposition factors can match the observed real data well, then the decomposition factors can be converted to original data structure to predict the missing entries.

In our algorithm, TT-WOPT is applied to real-valued tensor $\tensor{X}  \in\mathbb{R}^{I_1\times I_2\times\cdots \times I_N}$ with missing entries. The index of missing entries can be recorded by a weight tensor $\tensor{W}$ which is the same size as $\tensor{X}$. Every entry of $\tensor{W}$ meets:
\begin{equation}
\label{weight}
 w_{i_{1}i_{2}\cdots i_{N}}=
 \left\{
 \begin{aligned}
 &0 \qquad \text{if} \; x_{i_{1}i_{2}\cdots i_{N}} \;\text{is missing entry},\\
 &1 \qquad \text{if} \;  x_{i_{1}i_{2}\cdots i_{N}}\;\text{is observed entry}.
\end{aligned}
\right.
\end{equation}

In the optimization algorithm, the objective variables are the elements of all the core tensors. Define\, $\tensor{Y}=\tensor{W}\ast\tensor{X}$ and $\tensor{Z}=\tensor{W}\ast \ll \tensor{G}^{(1)},\tensor{G}^{(2)},\cdots,\tensor{G}^{(N)}\gg$ ($\ast$ is the Hadamard product, see \cite{jour6}), then the objective function can be written as:
\begin{equation}
\label{objective_function}
f(\tensor{G}^{(1)},\tensor{G}^{(2)},\cdots,\tensor{G}^{(N)})=\frac{1}{2} \left \|(\tensor{Y}-\tensor{Z}) \right \|^{2}.
\end{equation}

The relation between original tensor and core tensors can be deduced as the following equation \cite{book}:
\begin{equation}
\label{XG_relation}
\mat{X}_{(n)}=\mat{G}_{(2)}^{(n)}(\mat{G}_{(1)}^{>n} \otimes \mat{G}_{(n)}^{<n}),
\end{equation}
where for $n=1,...,N$, 

\begin{equation}
\mat{G}^{>n}=\ll \tensor{G}^{(n+1)},\tensor{G}^{(n+2)},\cdots,\tensor{G}^{(N)} \gg \in\mathbb{R}^{R_n\times I_{n+1}\times\cdots \times I_N},
\end{equation} 
\begin{equation}
\mat{G}^{<n}=\ll \tensor{G}^{(1)},\tensor{G}^{(2)},\cdots,\tensor{G}^{(n-1)} \gg \in\mathbb{R}^{I_1\times \cdots \times I_{n-1} \times R_{n-1}},
\end{equation}
where $\mat{G}^{>N}=\mat{G}^{<1}=E$ and $\otimes$ is the symbol of Kronecker products, also see \cite{jour6}.

For $n=1,...,N$, the partial derivatives of the objective function w.r.t. the $n$th core tensor $\tensor{G}^{(n)}$ can be inferred as follow:
\begin{equation}
\label{der}
\frac{\partial{f}}{\partial{\mat{G}_{(2)}^{(n)}}}=(\mat{Z}_{(n)}-\mat{Y}_{(n)})(\mat{G}_{(1)}^{>n} \otimes \mat{G}_{(n)}^{<n})^ \mathrm{ T }.
\end{equation}

After the objective function and the derivation of gradient are obtained, we can solve the optimization problem by any optimization algorithms based on gradient descent method \cite{jour12}. The optimization procedure of the algorithm is listed in $Alg. 1$.

\begin{table}[H]
\begin{center}
\begin{tabular}{l}
\hline
\textbf{Algorithm 1} Tensor-train Weighted Optimization (TT-WOPT)\\
\hline
\textbf{Input}: an $N$-way incomplete tensor $\tensor{X}$ and a weight tensor $\tensor{W}$.\\
\textbf{Initialization}: core tensors $\tensor{G}^{(1)},\tensor{G}^{(2)},\cdots,\tensor{G}^{(N)} $of tensor $\tensor{X}$.\\
1. Compute $\tensor{Y}=\tensor{W}\ast \tensor{X}$.\\
\textbf{For each optimization iteration,}\\
2. Compute $\tensor{Z}=\tensor{W}\ast \ll \tensor{G}^{(1)},\tensor{G}^{(2)},\cdots,\tensor{G}^{(N)}\gg$.\\
3. Compute objective function: $f=\frac{1}{2}\left \| \tensor{Y} \right \|^{2}-<\tensor{Y},\tensor{Z}>+\frac{1}{2}\left \| \tensor{Z} \right \|$.\\
4. Compute all $\frac{\partial{f}}{\partial{\mat{G}_{(2)}^{(n)}}}=(\mat{Z}_{(n)}-\mat{Y}_{(n)})(\mat{G}_{(1)}^{>n} \otimes \mat{G}_{(n)}^{<n})^ \mathrm{ T }$.\\
5. Use optimization algorithm to update $\tensor{G}^{(1)},\tensor{G}^{(2)},\cdots,\tensor{G}^{(N)} $. \\\textbf{Until} reach optimization stopping condition. \\
\textbf{Return} core tensors $\tensor{G}^{(1)},\tensor{G}^{(2)},\cdots,\tensor{G}^{(N)} $.\\
\hline
\end{tabular}
\end{center}
\end{table}

\section{Experiments}
In \cite{jour1} where the CP-WOPT method is proposed, only three-way data is tested. When it comes to high dimension data, the performance of CP-WOPT will fall. This is not because of the optimization method but the nature limit of CP decomposition. In our paper, we test our TT-WOPT on different orders of synthetic data. Then we test our algorithm on real world image data. We also compare the performance of TT-WOPT with several state-of-the-art methods. $\tensor{W}$ is created by randomly setting some percentage of entries to zero while the rest elements remain one.

\subsection{Simulation Data}
We consider to use synthetic data to validate the effectiveness of our algorithm. Till now, there is few relevant study about applying tensor-train decomposition to data completion, so we compare our algorithm to two other state-of-the-art methods--CP weighted optimization (CP-WOPT) \cite{jour1} and Fully Bayesian CP Factorization (FBCP) \cite{jour2}. We randomly initialize the factor matrices of a tensor with a specified CP rank, then we create the synthetic data by the factor matrices. For data evaluation index, we use relative square error (RSE) which is defined as $\scriptsize RSE=\sqrt{\left \| \tensor{X}-\hat{\tensor{X}} \right \|^2 / \left \| \tensor{X}\right \|^2 }$ where $\scriptsize\hat{\tensor{X}}$ is the tensor of full entries generated by core tensors or factor matrices. $Table$ \ref{sim results}. shows the simulation results of a three-way tensor and a seven-way tensor. The tensor sizes of synthetic data are $30\times 30 \times 30$ and $4 \times4 \times4 \times4 \times4 \times4 \times4 $, and the CP ranks are set to 10 in both cases.

Though we test the three algorithms on the data generated by CP model, our TT-WOPT algorithm shows good results. As we can see from Table \ref{sim results}., when we test on three-way tensor, TT-WOPT shows better fitting performance than CP-WOPT and FBCP at low data missing rates but a little weak at high missing rates. However, when we test on seven-way tensor, TT-WOPT outperforms the other two algorithms. In addition, we also find that the performance of TT-WOPT is sensitive to the setting of TT-ranks, different TT-ranks will lead to very different model accuracies. It should be noted that till now there is no good strategy to set TT-ranks and so in our experiments we set all TT-ranks the same value. This is an aspect that our algorithm needs to improve.  Furthermore, the initial values of core tensors also influence the performance of TT-WOPT.

\begin{table}[H]
\centering
\caption{Comparison of RSE of three different algorithms for two different data sizes with different missing rates of synthetic data. The algorithms are TT-WOPT, CP-WOPT, FBCP. The tensor ranks of each algorithm are set by experience (FBCP sets CP ranks automatically). The two data sizes are: $30\times 30 \times 30$ and $4 \times4 \times4 \times4 \times4 \times4 \times4 $. The three different missing rates are: 0\%, 50\% and 95\%.}
\label{sim results}
\resizebox{\textwidth}{!}{%
\begin{tabular}{c|c|c|c|c|c|c|c}
\hline
\hline
\multicolumn{2}{c|}{}  &\multicolumn{3}{c|}{three-way tensor} &\multicolumn{3}{c}{seven-way tensor} \\
\cline{1-8}
\multicolumn{2}{c|}{missing rate}   &0\% & 50\% &95\% &0\% &50\% &95\% \\
\hline
TT-WOPT&\makecell[cc]{TT-ranks \\ RSE}&\makecell[cc]{\{1,20,20,1\}\\ \textbf{2.64e-08} }&\makecell[cc]{\{1,20,20,1\}\\ \textbf{ 6.64e-05} }&\makecell[cc]{\{1,20,20,1\}\\1.06 }&\makecell[cc]{\{1,20,...,20,1\}\\ \textbf{7.22e-03} }&\makecell[cc]{\{1,20,...,20,1\}\\ \textbf{4.71e-03} }&\makecell[cc]{\{1,8,...,8,1\}\\0.744 }\\
\hline
CP-WOPT&\makecell[cc]{CP rank(manual) \\RSE}&\makecell[cc]{10\\ 5.34e-07}&\makecell[cc]{10 \\0.956 }&\makecell[cc]{10\\0.948 }&\makecell[cc]{10\\0.764}&\makecell[cc]{10\\0.957}&\makecell[cc]{10\\0.916 }\\
\hline
FBCP&\makecell[cc]{CP rank(auto)\\ RSE}&\makecell[cc]{10\\0.0581 }&\makecell[cc]{ 11\\0.0863 }&\makecell[cc]{6\\ \textbf{0.696}}&\makecell[cc]{3 \\ 0.542}&\makecell[cc]{7\\0.211 }&\makecell[cc]{4\\ \textbf{0.672} }\\
\hline
\hline
\end{tabular}}
\end{table}

\subsection{Image Data}
\begin{figure}[h]
\begin{center}
\includegraphics[width=0.65\linewidth]{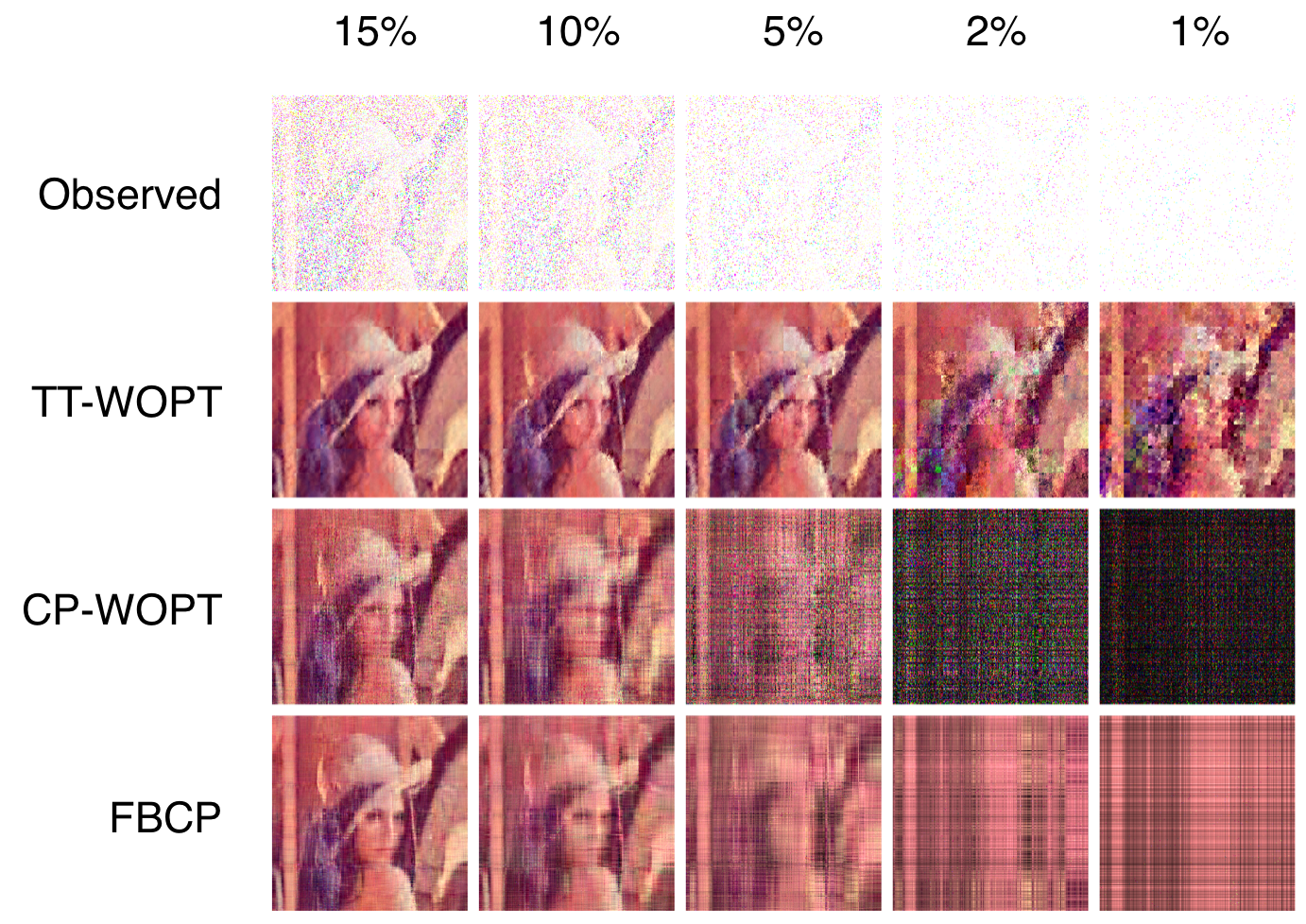}
\caption{Visualizing results of image inpainting performance of three different algorithms under five different missing rates of 85\%, 90\%, 95\%, 98\%  and 99\%. The values of missing entries of the image are changed from 0 to 255 in order to show the observed image clearly on the white paper.}
\label{default}
\end{center}
\end{figure}

In this section, we compare our algorithm with CP-WOPT and FBCP on image completion experiments. The size of every image data is $256\times 256 \times 3$. We use a set of images with missing rate from 85\% to 99\% to compare the performance of every algorithm. In this experiment, we do not set tensor ranks and tensor orders identically but use the best ranks to see the best possible result of every algorithm. For TT-WOPT, we first reshape  original data to a seventeen-way tensor of size $2 \times2 \times2 \times2 \times2 \times2 \times2 \times2 \times2 \times2 \times2 \times2 \times2 \times2 \times2 \times2 \times3 $ and permute the tensor according to the order of $\{1 \;9 \;2 \;10 \;3 \;11 \;4 \;12 \;5 \;13 \;6 \;14 \;7 \;15 \;8 \;16 \;17\}$. Then we reshape the tensor to a nine-way tensor of size $4 \times4 \times4 \times4 \times4 \times4 \times4 \times4 \times3$. This nine-way tensor is a better structure to describe the image data. The first-order of the nine-way tensor contains the data of a $2 \times 2$ pixel block of the image and the following orders of the tensor describe the expanding pixel blocks of the image. Furthermore, we set all TT-ranks to 16 according to our testing experience. For image evaluation index, we use PSNR (Peak Signal-to-noise Ratio) to measure the quality of reconstructed image data. $Table$ \ref{testing_results}. shows the testing results of one image. $Fig.$\ref{default}. visualizes the image inpainting results.

The experiment result shows that our TT-WOPT algorithm outperforms other algorithms for image data completion. Particularly, when the missing rate reaches 98\% and 99\%, our algorithm can recover the image successfully while other algorithms totally fail. The RSE and the PSNR values of TT-WOPT are always better than CP-WOPT and FBCP. In addition, the image visual quality of our method is always the best. 

\begin{table}[h]
\centering
\caption{Comparison of the inpainting performance (RSE and PSNR) of three algorithms under five different missing rates: 85\%, 90\%, 95\%, 98\% and 99\% of a testing image. }
\label{testing_results}
\begin{tabular}{c|c|c|c|c|c|c}
\hline\hline
 \multicolumn{2}{c|}{missing rate}& 85\% & 90\% & 95\% & 98\% & 99\% \\
 \hline
TT-WOPT&\makecell[cc]{RSE \\ PSNR}   &     \makecell[cc]{\textbf{0.1233} \\ \textbf{23.4877}}    &    \makecell[cc]{\textbf{0.1297} \\ \textbf{22.6076}}     &
\makecell[cc]{\textbf{0.1416}\\ \textbf{21.5282}}     &           \makecell[cc]{\textbf{0.2202}\\ \textbf{18.9396}}   &\makecell[cc]{\textbf{0.2638}\\ \textbf{17.0029}}\\
\hline
CP-WOPT&\makecell[cc]{RSE \\ PSNR}&\makecell[cc]{0.1891\\ 18.8578}&\makecell[cc]{0.3169 \\ 18.0389}&\makecell[cc]{0.5348 \\ 12.5649}&\makecell[cc]{1.0918 \\7.8015} &   \makecell[cc]{1.1309\\6.4971}\\
\hline
FBCP&\makecell[cc]{RSE \\ PSNR}&\makecell[cc]{0.1440\\ 22.2853}&\makecell[cc]{0.1867 \\ 19.9410}&\makecell[cc]{0.2432\\ 17.5166}&\makecell[cc]{0.3052 \\ 15.4784}  & \makecell[cc]{0.3372\\ 14.5841}\\
\hline\hline
\end{tabular}
\end{table}

\section{Conclusion}
In this paper, we first elaborate the basis of tensor and the tensor-train decomposition method. Then we use a gradient-based first-order optimization method to find the factors of the tensor-train decomposition when tensor has missing entries and propose the TT-WOPT algorithm. This algorithm can solve the tensor completion problem of high dimension tensor. From the simulation and image experiments, we can see our algorithm outperforms the other state-of-the-art methods in many situations especially when the missing rate of data is extremely high. Our study also proves that high order tensorization of data is an effective and efficient method to represent data. Furthermore, it should be noted that the accuracy of TT model is sensitive to the selection of TT-ranks. Hence, we will study on how to choose TT-ranks automatically in our future work.

\section*{Acknowledgement} 
This work is supported by JSPS KAKENHI (Grant No. 17K00326) and KAKENHI (Grant No. 15H04002).

\end{document}